\documentclass{ifacconf}

\usepackage{amsmath,amssymb}
\usepackage{graphicx}   
\usepackage[square,numbers]{natbib}        
\usepackage{algorithm}
\usepackage{algorithmic}

\pdfoutput=1
\begin{document}
\begin{frontmatter}

\title{On Optimal Control of Hybrid Dynamical Systems using Complementarity Constraints} 

\thanks[footnoteinfo]{Work at Los Alamos National Laboratory is done under the auspices of the National Nuclear Security Administration under U.S. D.O.E. Contract No. 89233218CNA000001.}

\author[First]{Saif R. Kazi} 
\author[Second]{Kexin Wang}
\author[Third]{Lorenz Biegler}

\address[First]{T-5 Applied Mathematics and Plasma Physics, Los Alamos National Laboratory, Los Alamos, NM 87545, USA (e-mail: skazi@lanl.gov)}
\address[Second]{College of Control Science and Engineering, Zhejiang University, Hangzhou, Zhejiang 310027, China (e-mail: kxwang@zju.edu.cn)}
\address[Third]{Department of Chemical Engineering, Carnegie Mellon University, Pittsburgh, PA 15213, USA (e-mail: lb01@andrew.cmu.edu)}

\begin{abstract}                
Optimal control for switch-based dynamical systems is a challenging problem in the process control literature. In this study, we model these systems as hybrid dynamical systems with finite number of unknown switching points and reformulate them using non-smooth and non-convex complementarity constraints as a mathematical program with complementarity constraints (MPCC). We utilize a moving finite element based strategy to discretize the differential equation system to accurately locate the unknown switching points at the finite element boundary and achieve high-order accuracy at intermediate non-collocation points. We propose a globalization approach to solve the discretized MPCC problem using a mixed NLP/MILP-based strategy to converge to a non-spurious first-order optimal solution. The method is tested on three dynamic optimization examples, 
including a gas-liquid tank model and an optimal control problem with a sliding mode solution. 



\end{abstract}


\end{frontmatter}

\section{Introduction}

Hybrid Dynamical Systems have become a popular method to model systems that are characterized by a mixed continuous and discrete solution to the physical state variables \cite{van2007introduction}. Common examples with process control applications include robot dynamics \cite{patel2019contact,raghunathan2022pyrobocop,wei2008optimal}, chemical processes \cite{Baumrucker2009mpec,raghunathan2004MPEC,moudgalya2001class-II} and other areas \cite{biak2013some,nurkanovic2022nosnoc}. The problem is computationally difficult to solve as it relates to optimization of a dynamical system with discrete switch variables. 

Multiple studies have proposed different ways to solve the optimal control of hybrid dynamical systems, including maximum principle approach \cite{Riedinger1999maximum,riedinger2003optimal,cassandras2001optimal}, Bellman inequality \cite{hedlund1999optimal}, mixed integer programming (MIP) \cite{bemporad1999control}, dynamic programming \cite{xu2004optimal} and direct transcription \cite{buss2002nonlinear,bengea2005optimal}. The last approach relates to numerical scheme where the dynamical differential equations are reformulated and discretized into a set of algebraic constraints. This approach leads to a non-smooth optimization problem which is either solved as a large-scale MIP or a mathematical program with complementarity constraints (MPCC) \cite{stewart1990high,stewart2010optimal}.  

The modeling of switching points and discretization for hybrid dynamical systems is complex and non-trivial. The location and number of switching points is usually unknown apriori. Also, higher order discretization schemes fail to achieve high accuracy solutions, namely, discretization error decreases slowly as a function of number of grid points or step size. This happens due to the non-smooth nature of the MPCC formulation and requires sophisticated discretization and solution techniques as compared to standard nonlinear programming (NLP) solvers. The formulation of hybrid dynamical systems and its reformulation to an MPCC are described in the rest of this section.



\subsection{Hybrid Dynamical Systems}

A typical Hybrid Dynamical System can be represented as ordinary differential equations 
with discontinuous right-hand side as:
\begin{equation}\label{eq:hybriddynamics}
    \dot{x} = f_j (x(t),u(t)), \text{ if } x(t) \in R_j \subset \mathbb{R}^{n_x}, \enskip j \in \mathcal{F} = \{1, \dots, n_f\},
\end{equation}

where $R_j$ denote the regions (disjoint sets) in the state space and
$f_j(.): \mathbb{R}^{n_x+n_u} \to \mathbb{R}^{n_x}$ are corresponding smooth functions representing the state
dynamics in them; $n_f$ is the number of piecewise functions and
$u(t)\in U \subseteq \mathbb{R}^{n_u}$ denotes the external control inputs in the system.
\subsection{Filippov Reformulation and Optimization}

Hybrid Dynamical Systems were extensively studied by
\citet{filippov1960differential}, who proposed a formulation using the
convex combination to transform \eqref{eq:hybriddynamics} as
\begin{equation}\label{eq:Filippov}
    \dot{x} = \sum_{j \in \mathcal{F}} \nu_j f_j (x,u), \enskip  \sum_{j \in \mathcal{F}} \nu_j = 1, \nu_j \geq 0, \nu_j = 0 \text{ if } x(t) \notin R_j,
\end{equation}
where $\nu_j$ are also called \textit{convexification variables} or \textit{Filippov
multipliers}. The convexificaton allows to write the piecewise
formulation \eqref{eq:hybriddynamics} with disjoint sets into a
continuous algebraic form \eqref{eq:Filippov} where the system
dynamics at the boundary region $\{ x \in R_{j'} \cap R_j \}$ is well
defined, i.e., $\dot x = \nu_{j'} f_{j'} + \nu_j f_j$, $ \nu_{j'} +
\nu_j = 1$, which is more suitable for continuous optimization
solvers. \citet{stewart1990high} presented an extension with
\textit{indicator} functions ($g_j(x)$) for the disjoint sets ($R_j$):
\begin{equation}
    R_{j'} = \{ x \in \mathbb{R}^{n_x} | g_{j'}(x) < \min_{j \in \mathcal{F}, \ne j'} g_j(x) \}. 
\end{equation}

Under mild assumptions such that $g(.)$ and $\nabla g(.)$ are Lipschitz continuous and sufficiently smooth, the Filipov system can be reformulated as:
\begin{subequations}
    \begin{equation}
        \dot x = \sum_{j \in \mathcal{F}} \nu_j f_j (x,u) ,
    \end{equation}
    where the \textit{Filippov multipliers} are algebraic variables determined by the following parametric linear program (LP):
    \begin{equation}\label{eq:LP}
        \text{LP}(x): \enskip \nu(x) = \arg\min g(x)^T \nu \text{ s.t. } \sum_{j \in \mathcal{F}} \nu_j = 1, \nu_j \geq 0 .
    \end{equation} 
\end{subequations}
Since the parametric LP solution vector ($\nu(x)$) is dependent on the state variable ($x$), which itself is a primal variable in the optimal control problem. The linear program \eqref{eq:LP} needs to be reformulated using its KKT conditions for the optimal control problem to be represented as a single-level optimization problem. Thus, the hybrid dynamic constraints in the problem can be rewritten into a dynamic complementarity system (DCS) as follows:
\begin{subequations}\label{eq:DCS}
    \begin{align}
        \dot x = \sum_{j \in \mathcal{F}} \nu_j f_j (x,u),\\
        g(x) - \lambda(t) - \mu(t) e = 0,\\
        \sum_{j \in \mathcal{F}} \nu_j(t) = 1,\\
        0 \leq \nu(t) \perp \lambda(t) \geq 0 ,
    \end{align}
    where $g(x) = [g_1(x),\dots, g_{n_f}(x)]^T, e = [1,\dots,1]^T$ denote
the indicator functions and the unit vector, whereas $\nu(t)$ denotes the 
convexification variables, and $\lambda(t) \in \mathbb{R}^{n_f}, \mu(t) \in \mathbb{R}$ denote the
KKT multipliers; all are functions of $t$. 
Complementarity constraints are written as $0 \leq x \perp y \geq 0$, 
where $\perp$ denotes the complementarity operator. 

\end{subequations}

The complementarity constraint (\ref{eq:DCS}d) between the inequality multipliers $\lambda$ and Filippov variables $\nu$ is equivalent to:
\begin{equation}\label{eq:MPCC-equivalence}
    0 \leq \nu \perp \lambda \geq 0  \Longleftrightarrow \nu, \lambda \geq 0, \; \nu^T \lambda = 0 .
\end{equation}
Since both vectors are non-negative, the relation is equivalent to
$\nu_j \lambda_j = 0$ or $\min(\lambda_j,\nu_j) = 0 \enskip \forall j
\in \mathcal{F}$.

The optimal control problem for hybrid dynamic systems is described as an infinite dimensional optimization problem with objective function consisting of terminal cost function $\phi$ and stage cost $\Psi$ constrained to differential algebraic equations Eq.\ref{eq:dynamic_infinite_opt}b-\ref{eq:dynamic_infinite_opt}d and with the complementarity constraints Eq.\ref{eq:dynamic_infinite_opt}e as shown here:
\begin{subequations}\label{eq:dynamic_infinite_opt}
    \begin{align}
        \min_{u} \enskip & \phi (x(t_f)) + \int^{t_f}_0 \Psi (x,u) dt \\
        \text{s.t.} \enskip & \dot x = \sum_{j \in \mathcal{F}} \nu_j f_j (x,u), \enskip  x(0) = x_0,\\
        & g(x) - \lambda - \mu e = 0,\\
        & \sum_{j \in \mathcal{F}} \nu_j = 1,\\
        & 0 \leq \nu \perp \lambda \geq 0, \qquad t \in [0,t_f] .
    \end{align}
\end{subequations}

\subsection{Challenges}

There are two major challenges in solving the dynamic optimization with complementarity constraints as formulated in \eqref{eq:dynamic_infinite_opt}. 

\subsubsection{Switching Point Location.}

The switching point(s) for hybrid dynamical systems are unknown a priori and are functions of the state variable ($g_j(x)$). Solving the dynamic optimization problem (\ref{eq:dynamic_infinite_opt}) with a higher order numerical discretization scheme requires non-uniform discretization to accurately locate the switching points. The moving finite element strategy to precisely detect the switching points is discussed in Section 2.

\subsubsection{Complementarity Constraints.}

The complementarity constraint in (\ref{eq:dynamic_infinite_opt}e) is
nonlinear, non-convex and non-differentiable. This makes the dynamic
optimization problem difficult to solve with gradient-based
optimization solvers. Necessary conditions such as constraint
qualifications, i.e. LICQ and MFCQ are violated at each feasible point.
A detailed discussion on MPCCs is presented in Section 3. 

\subsection{Outline}

The rest of the paper is outlined as follows. 
\begin{itemize}
\item
Section 2 describes the discretization and moving finite element strategy to accurately detect the switching points in the optimal solution.  This section proposes a discretization method which shares the high accuracy of the moving finite element strategy while consuming the additional degrees of freedom in step size without increasing numerical difficulty.
\item
Section 3 discusses optimization with complementarity constraints. This section develops an algorithm that seeks a B-stationary solution of MPCCs. It can be expected to perform efficiently, because every trial point is feasible to the MPCC and thus allows identification of the active sets;
descent directions are determined based on the bi-active instead of the whole set of the complementarity constraints;
objective improvement is guaranteed at every trial point in order to enforce convergence of the B-stationarity exploration.
\item
In Section 4, the proposed algorithm is implemented on three case studies to demonstrate the effectiveness of the approach. 
\item
Conclusions and recommendations for future work are discussed in Section 5.
\end{itemize}

\section{Non-smooth Dynamic Optimization and Formulations}

The non-smooth dynamic optimization with complementarity constraints
presented in (\ref{eq:dynamic_infinite_opt}) can be solved as a
differential algebraic equation (DAE) optimization with the
complementarity constraints reformulated as equality and inequality
constraints as shown in \eqref{eq:MPCC-equivalence}. DAE optimization
problems can be approached in two primary ways: 1) \textit{optimize -
then - discretize} and 2) \textit{discretize - then - optimize}. The
first method formulates the optimization problem in continuous form and solves its
first-order optimality conditions (Euler-Lagrange equations) as a coupled forward-backward adjoint boundary value
problem (BVP). While this approach 
is derived from the calculus of variations in the continuous time domain, it may be challenging to
implement for general dynamical constrained optimization problems.

\subsection{Discretization Methods}

The second approach \textit{discretize - then - optimize} is a more common and generalizable method to solve dynamic optimization problems. The differential equation system is discretized using numerical discretization schemes such as Runge-Kutta (RK) or orthogonal collocation matrix (OCM) over finite elements (FE). This results in a discretized system of nonlinear algebraic equations which is then solved as a large-scale nonlinear programming (NLP) problem. The numerical accuracy of the NLP solution using this approach can be increased using higher order discretization schemes such as RK4 or 3 point Radau collocation method. 

\begin{subequations}\label{eq:standard IRK}

    For brevity, we denote $F = [f_1,\dots, f_{n_f}]$ and rewrite
   (\ref{eq:dynamic_infinite_opt}b) and (\ref{eq:dynamic_infinite_opt}d) as $\dot x = F(x,u)^T \nu, $ and
    $e^T \nu = 1$. \\ \\
    For $l = 1,\dots,N, \enskip k = 1,\dots,K,$
    \begin{align}
        x_{l+1,0} = x_{l,0} + h \sum^K_{k=1} b_k F_{l,k} \nu_{l,k}, \\
        x_{l,k} = x_{l,0} + h \sum^K_{k=1} a_{l,k} F_{l,k} \nu_{l,k},\\
        g(x_{l,k}) - \lambda_{l,k} - \mu_{l,k}e = 0,\\
        e^T \nu_{l,k} = 1,\enskip x_{0,0} = x_0,\\
        0 \leq \nu_{l,k} \perp \lambda_{l,k} \geq 0 ,
    \end{align}
    where $N$ is the number of finite elements and $K$ is the number of
collocation points or stages; $a_{l,k}$ and $b_k$ are constant
parameters depending on the discretization scheme; $h$ refers to the
finite element step size.
\end{subequations}

\subsection{Moving Finite Element Strategy}

For typical DAE optimization problems, the solution accuracy is
primarily a function of step size ($h$) and discretization order ($K$).
On the other hand, for DAE optimization with complementarity constraints, the
higher solution accuracy is not guaranteed for uniform grid
discretization, i.e., constant step size ($h$). The solution may have
lower numerical accuracy and higher discretization error, if the
switching point between the complementarity variables ($\nu_{l,k}$,
$\lambda_{l,k}$) is not at a finite element boundary. 

\begin{figure}[h]
    \centering
    \includegraphics[width=1.0\linewidth]{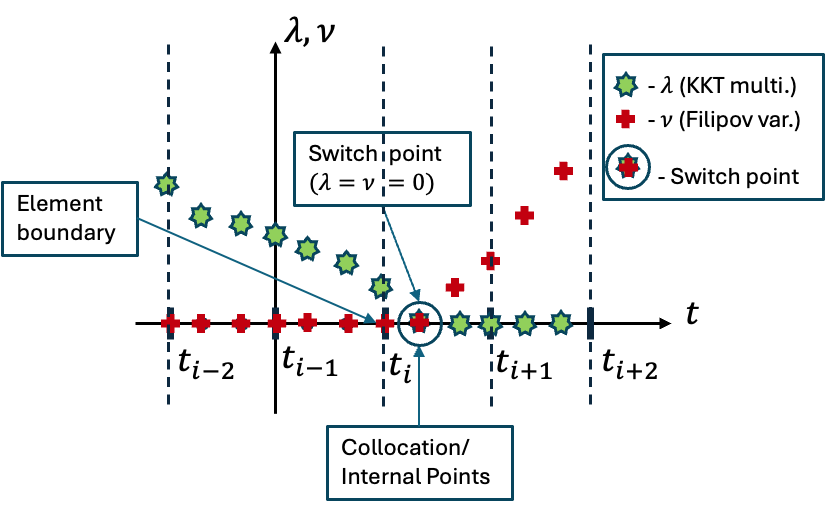}
    \caption{Illustration shows that the non-differentiable switching point for uniformly discretized solution may exist inside the finite element, thus violating the Taylor series assumption.}
    \label{fig:smoothness violation}
\end{figure}

As seen in Fig. \ref{fig:smoothness violation}, the complementarity condition between the variables is satisfied at each collocation point, but the switch between the complementarity variables, i.e., the point where either $\nu$ or $\lambda$ switch from being strictly positive to zero or vice-versa, is inside the finite element $[t_l, t_{l+1}]$. As the solution at the switch point is non-differentiable, the assumption of discretization methods, based on Taylor series of continuously smooth solutions, is violated inside the finite element. This increases the discretization error and reduces the accuracy order of the numerical discretization scheme with respect to step size.  

Baumrucker and Biegler \cite{Baumrucker2009mpec} presented a moving finite element (MFE) strategy where the discretization and step sizes $h_l$ are non-uniform and declared as variables in the optimization problem with $h_l \in [\underline{h},\bar h]$.  

\begin{subequations}\label{eq:cross-complementarity}
    For $l = 0,\dots,N-1, \enskip k = 1,\dots,K,$
    \begin{align}
        x_{l+1,0} = x_{l,0} + h_l \sum^K_{k=1} b_k F_{l,k} \nu_{l,k}, \\
        x_{l,k} = x_{l,0} + h_l \sum^K_{k=1} a_{l,k} F_{l,k} \nu_{l,k},\\
        \sum^{N-1}_{l=0}h_l = t_f,  \underline{h} \leq h_l \leq \bar{h},\\
        g(x_{l,k}) - \lambda_{l,k} - \mu_{l,k}e = 0,\\
        e^T \nu_{l,k} = 1,\enskip x_{0,0} = x_0,\\
        0 \leq \nu_{l,k} \perp \sum^K_{k=1} \lambda_{l,k} \geq 0 .
    \end{align}
    The bounds on the step size ensure that the element boundaries
don't coincide with each other and the problem is well-conditioned.
\end{subequations}

To ensure that the switch point is located at the boundary of the
finite element and the solution remains continuously differentiable
inside the finite element, \citet{Baumrucker2009mpec} reformulated the
complementarity constraints in (\ref{eq:standard IRK}e) into
cross-complementarity constraints (\ref{eq:cross-complementarity}f).
This reformulation complements the sum of the variables ($\sum_k
\lambda_{l,k}$) with each complementing variable $\nu_{l,k}$ inside
the element. This ensures that $\nu_{l,k} > 0$ is allowed only if $\sum_k     
\lambda_{l,k} = 0$, and the complementarity variables do not
switch inside the finite element, as shown in Fig.\ref{fig:moving
finite element}.

\begin{figure}[h]
    \centering
    \includegraphics[width=1.0\linewidth]{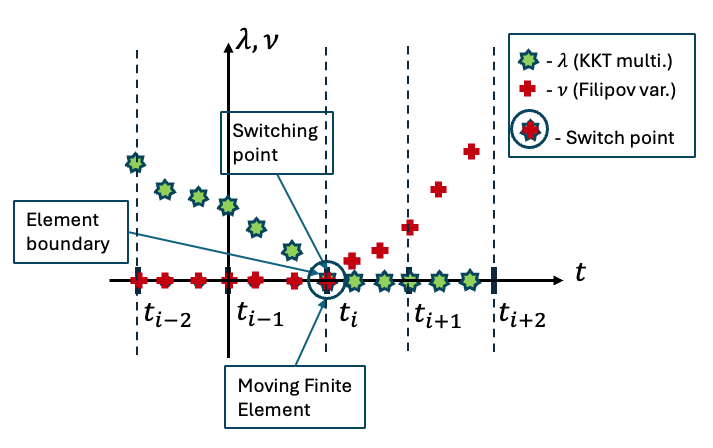}
    \caption{Illustration shows how the moving finite element strategy ensures that the non-differentiable switching point coincides with the boundary point of the finite element, thus satisfying the Taylor series assumption.}
    \label{fig:moving finite element}
\end{figure}

\subsection{Nurkanovic Formulation}

Although the moving finite element reformulation ensures that the
switching points in the solution coincide with the element boundaries,
the additional variables $h_l$ and only
one additional equality constraint (\ref{eq:cross-complementarity}c)
in the formulation increases the number of free variables, which may
result in non-unique values for discretization variables $h_l$. 

\citet{nurkanovic2022finite} proposed additional constraints to ensure
that the 
discretization values $h_l$ are unique. The proposed
step-equilibration approach is based on the principle that only the
finite element(s) with the switch point(s) have non-uniform
discretization $[t_{l-1}, t_{l}]$. The approach uses an indicator
variable ($\eta$) for the switch point(s), which is determined by
calculating and multiplying the sum of the complementarity variables
in two consecutive finite elements. 

\begin{subequations}\label{eq:step-equilibration}
The auxiliary variables for the sum of the complementarity variables at each finite element are defined as:
    \begin{align}
        \hat{\lambda_l} = \sum^K_{k=1} \lambda_{l,k}, \quad \hat{\nu_l} = \sum^K_{k=1} \nu_{l,k} .
    \end{align}
    Then, the Hadamard product of the forward and backward sum of the complementarity variables determine if they have switched from positive to zero (or vice-versa).
    \begin{align}
        \pi^{\lambda}_l = \hat{\lambda}_{l-1} \odot \hat{\lambda_l}, \quad \pi^{\nu}_l = \hat{\nu}_{l-1} \odot \hat{\nu_l} .
    \end{align}
    (Here, $\odot$ represents pointwise or elementwise product of vectors.)\\
    Since at least one of the vectors $\pi^{\lambda}_l$ or $\pi^{\nu}_l$ is zero at each element, and they are exactly equal to zero at the element corresponding to the switching point, the sum of the two vectors is a good candidate for the indicator function
    \begin{align}
        \tau_l &= \pi^{\lambda}_l + \pi^{\nu}_l, \; \eta_l = \prod^{n_f}_{j=1} \tau_{l,j}.
    \end{align}
    Since the indicator variable $\eta_l$ is non-negative and only zero at the switching element, the relation between step size and indicator variable can be represented by the following complementarity constraints.
    
    For $l = 1,\dots,N-1,$
    \begin{align}
        0 \leq (\Delta  h_{l}^+ + \Delta  h_{l}^-) \perp \eta_l \geq 0,
    \end{align}
where $h_{l-1} - h_{l} = \Delta  h_{l}^+ - \Delta  h_{l}^-$, $\Delta  h_{l}^+, \Delta  h_{l}^- \geq 0$,
\end{subequations}
The finite element with switch detection (FESD) algorithm was implemented as a package NOSNOC in \cite{nurkanovic2022nosnoc}.

As mentioned in \eqref{eq:step-equilibration}, the Nurkanovic formulation augments an additional  [$2N \{\hat \lambda, \hat \nu\} + (2N-2) \{\pi^{\lambda},\pi^{\nu}\} + (N-1) \{\tau\} + (N-1) \{\eta\}$] variables and [$2N \{(\ref{eq:step-equilibration}{\rm a})\} + (2N-2) \{(\ref{eq:step-equilibration}{\rm b}) \} + 2(N-1) \{(\ref{eq:step-equilibration}{\rm c}) \} +  (N-1) \{(\ref{eq:step-equilibration}{\rm d})\}$] constraints for each complementarity constraint. This effectively decreases the degrees of freedom by $N-1$ to that of the original problem and avoids non-unique solutions for the step size variables $h_i$. 

\subsection{Proposed Two Stage Step Equilibration Approach}

Although the Nurkanovic formulation makes the problem consistent with respect to the degrees of freedom and ensures uniform grid discretization away from the switch point(s), the formulation may be numerically unstable (i.e. the auxiliary variables have small values, the Jacobian or Hessian has large condition number) and increases the size of the problem, making it difficult to implement on larger optimal control problems. 

Prompted by the Nurkanovic \cite{nurkanovic2022finite} formulation, we propose a modification of the approach in \cite{Baumrucker2009mpec}, in order to keep the degrees of the problem consistent. Our proposed two step approach is described below:

\textbf{First Stage:}
In the first stage, we apply the moving finite element (MFE) formulation in \cite{Baumrucker2009mpec} with cross complementarities (\ref{eq:cross-complementarity}f) and solve it to optimal point. The MFE method ensures that the switching point coincides with the boundary of the finite elements. The switching point(s) ($t_s$) are located by observing the solution values of $\lambda$ and $\nu$ at the boundary of each finite element.

Define the set of finite elements which have the switching point at the end (i.e. right) as:
\begin{subequations}\label{eq:proposed-formulation}
\begin{align}\label{eq:switch-point}
\begin{split}
   & \chi_s = \{l \enskip | \enskip [t_{l-1}, t_l], \enskip \lambda_{l,K} = \nu_{l,K} = 0,\\
   \enskip & \lambda_{l-1,K} + \nu_{l-1,K} > 0 \enskip \text{or} \enskip \lambda_{l+1,K} + \nu_{l+1,K} > 0 \} .
\end{split}
\end{align}

In Eq.\eqref{eq:switch-point}, we compare the switching variables ($\lambda_{l,K}$ and $\nu_{l,K}$) at the boundary of consecutive finite elements. Switching point is identified if and only if both the switching variables are equal to 0 ($\leq \delta_{tol}$) at a finite element boundary, but both are not simultaneously equal to 0 at the previous or next finite element. This ensures that we only count the starting and end points of switching for sliding mode solutions.  

\textbf{Second Stage:}
In the next stage, we add additional constraints to the MFE formulation which forces the finite elements to be equally spaced away from the switching points.

For finite elements with no switching points at their boundary, we enforce the uniform discretization using the following constraint: 
\begin{align}\label{eq:uniform_step}
    h_l - h_{l+1} = 0, \quad \forall l \in \{1,\dots,N-1\} \setminus \chi_s .
\end{align}

For finite elements with switching points, we add constraints to force the switching to happen at the boundary of finite elements found in the first stage.
\begin{align}\label{eq:non-uniform_step}
    \lambda_{l,K} + \nu_{l,K} = 0, \quad \forall l \in \chi_s .
\end{align}
\end{subequations}

The proposed formulation adds the necessary $N-1$ linear constraints without the additional variables used in Nurkanovic formulation. The resultant optimal control problem is thus more stable and numerically well-conditioned for nonlinear optimization solvers. We solve the resultant problem with Eq.(\ref{eq:cross-complementarity}), Eq.(\ref{eq:uniform_step}) and Eq.(\ref{eq:non-uniform_step}). 

The main assumption in our approach is that the location of the switching points in the optimal solution is independent of the step-size variables and the formulations. Thus, the switching points in the MFE formulation would be the same as in the Nurkanovic formulation. To obtain unique solution and avoid complex constraints in Nurkanovic formulation, we implement uniform discretization between the start time, switch points and the final time using \eqref{eq:proposed-formulation} instead.  

\section{Solution methods for MPCCs}

To develop the solution strategy for the MPCC derived in the previous section, we discretize and rewrite the optimal control problem \eqref{eq:dynamic_infinite_opt}
in the more general form:
\begin{subequations}\label{eq:general-MPCC}
    \begin{align}
        \min \enskip & \varphi(x)\\
        \text{s.t.} \enskip & g_I(x) \leq 0; \; g_E(x) = 0,\\
        & 0 \leq G(x) \perp H(x) \geq 0 ,
    \end{align}
\end{subequations}
where $g_I$ and $g_E$ represent the discretized equality and inequality constraints, and the complementarity constraints (\ref{eq:general-MPCC}c) represent the cross-complementarity constraints (\ref{eq:cross-complementarity}f). 

\subsection{Stationarities and First-Order Optimality}
The following index sets are defined at every feasible point $\bar x$ of the MPCC \eqref{eq:general-MPCC}:
\begin{subequations}
    \begin{align}
        I_{g_I}(\bar x) &= \{i| \enskip g_{I_i}(\bar x) = 0 \}, \\
        I_G(\bar x) &= \{i| \enskip G_i(\bar x) = 0, H_i (\bar x) > 0 \}, \\
        I_H(\bar x) &= \{i| \enskip G_i (\bar x) > 0, H_i(\bar x) = 0\}, \\
        I_{GH}(\bar x) &= \{i| \enskip G_i(\bar x) = 0, H_i(\bar x) = 0 \} .
    \end{align}
\end{subequations}
Note that the sets $I_{G}(\bar x)$ and $I_{H}(\bar x)$ are disjoint, and $I_{GH}(\bar x)$ is also called the bi-active set at point $\bar x$, 

Based on these sets, stationarity of $\bar x$ can be specified. In particular, $\bar x$ is W-stationary (weakly stationary) if there exist multipliers $\lambda^{g_I} \geq 0,\lambda^{g_E},\sigma^G,$ and $\sigma^H$, such that $\lambda^{g_I}_i g_{I_i}(\bar x) =0$ for $i = 1 \dots |g_I|$, $\sigma^G_i G_i(\bar x) =0$ and
$\sigma^H_i H_i(\bar x) =0$ for $i = 1 \dots n_c$, and
\begin{equation*}
\begin{aligned}
    0 &= \nabla \varphi(\bar x) 
       + \nabla g_I(\bar x) \lambda^{g_I}
       + \nabla g_E(\bar z) \lambda^{g_E} \\
      &- \sum_{l \in I_G(\bar x) \cup I_{GH}(\bar x)} 
       \sigma_i^G \nabla G_i(\bar x)
       - \sum_{i \in I_H(\bar x) \cup I_{GH}(\bar x)} \sigma_i^H \nabla H_i(\bar x) .
\end{aligned}
\end{equation*}
Moreover, a weakly stationary point $\bar x$ is also
\begin{itemize}
    \item
        A-Stationary, if either $\sigma_i^G \geq 0$ or $\sigma_i^H \geq 0$ for all $i \in I_{GH}(\bar x)$;
    \item
        C-stationary, if $\sigma_i^G \sigma_i^H \geq 0$ for all $i \in I_{GH}(\bar x)$;
    \item
        M-stationary, if either $\sigma_i^G,\sigma_i^H \geq 0$ or $\sigma_i^G \sigma_i^H = 0$ for all $i \in I_{GH}(\bar x)$;
    \item 
        S-stationary (strongly stationary), if $\sigma_i^G,\sigma_i^H \geq 0$ for all $i \in I_{GH}(\bar x)$.
\end{itemize}
Not all of these stationarity concepts are immediately helpful in recognizing a local minimum of an MPCC.
At first sight, W-, A-, C-, M- stationarities may not exclude the existence of first-order descent directions at $\bar z$, and therefore may not be first-order conditions for $\bar z$ to be a local minimum.
On the other hand, S-stationarity guarantees that no first-order descent directions exist. But considering only S-stationarity at $\bar z$, assumes that $\bar z$ is also first-order optimal even if the bi-active complementarity constraints are replaced by $G_i(x), H_i(x) \geq 0$ for all $i \in I_{GH}(\bar z)$. This assumption is clearly stronger than being first-order optimal for the MPCC.
In fact, S-stationarity is a first-order necessary condition, only if the constaint set satisfies certain regularity conditions at $\bar z$ (see, for example, \cite{scheel2000mathematical, bertsekas-ozdaglar}). 
For more general cases, we need the concept of B-stationarity to confirm first-order optimality.

A feasible point $\bar x$ of the MPCC (\ref{eq:general-MPCC}) is B-stationary, if it satisfies
\begin{equation} \label{bstat}
    \nabla \varphi(\bar x)^Td \geq 0, \quad
    \forall d \in \mathcal T(\bar x),
\end{equation}
where $\mathcal T(\bar x)$ is the tangent cone of the constraint set at $\bar x$.
This definition basically says that at a B-stationary point $\bar x$, there are no feasible directions (that is, directions in the tangent cone) that can decrease the objective $\varphi(\bar x)$.
Therefore, B-stationarity is a precise necessary condition for first-order optimality. 

As a property pertaining to every local minimum of an MPCC, B-stationarity can take various concrete forms.
It is equivalent to S-stationarity if the MPCC linear independence constraint qualification (MPCC-LICQ) holds at $\bar x$ \cite{scheel2000mathematical} (i.e., the equality constraints and active inequality constraints in (\ref{eq:general-MPCC}) are linearly independent).
Otherwise, it can be strictly weaker than S-stationarity. For instance, Examples 2.1, 2.3, and
5.2 in \cite{guo-lin-ye} and ex9.2.2 from \cite{leyffer2000} have local minima that are no better than M-stationary, and Example 2.4 in \cite{guo-lin-ye} has an unique minimum that is no better than C-stationary.
The complex relationship between B-stationarity and other stationarities makes direct determination of B-stationary points difficult. Therefore, we require an alternate approach to obtain a B-stationary point, which is developed in the remainder of this section. 

\subsection{NLP-based Strategies and Convergence Properties}

The following are some representative NLP-based solution strategies
for MPCCs. These strategies are designed not to deal with the
complementarity structure explicitly. As a consequence,
we can see from their convergence properties the limitations in guaranteeing first-order optimality.

\subsubsection{Regularization.}
The typical regularization scheme proposed by \cite{scholtes2001} reformulates (\ref{eq:general-MPCC}) into the form
\begin{equation*}
\begin{aligned}
    {\rm REG}(\epsilon): \; 
    \min \; & \varphi(x) &{\rm multipliers}\\
    {\rm s.t.} \; &g_I(x) \leq 0, &v^{g_I} \\
    ~ &g_E(x) = 0, &v^{g_E} \\
    ~ &G(x) \geq 0, &v^G \\
    ~ &H(x) \geq 0, &v^H \\
    ~ &G_i(x)H_i(x) \leq \epsilon, \; i= 1 \dots n_c. &v^{REG}_i
\end{aligned}
\end{equation*}
Solving a sequence of problems REG($\epsilon^k$) with the positive scalars $\epsilon^k \to 0$, yields a sequence of stationary points $x^k$ of REG($\epsilon^k$), which tends to an accumulation point $\bar x$.
Using the Lagrange function defined by: 
\begin{equation*}
\begin{aligned}
&\mathcal L = \varphi(x) 
+ (v^{g_I})^T g_I(x) 
+ (v^{g_E})^T g_E(x) \\
&- (v^{G})^T G(x)- (v^{H})^T H(x)+ \sum_{i=1}^{n_c} v_i^{REG}(G_i(x)H_i(x) - \epsilon) .
\end{aligned}
\end{equation*}
B-stationarity of $\bar x$ has been established under MPCC-LICQ at
$\bar x$ and the second-order necessary conditions at every $x^k$,
together with the assumption of upper level strict complementarity,
namely, $v^G_i v^H_i \neq 0$ for all $i \in I_{GH}(\bar x)$.


\subsubsection{NCP Function.}
An NCP function represents a complementarity constraint with a suitable nonlinear and usually nondifferentiable equation. Here we consider the following NCP function with a smoothing factor $\epsilon >0$:
\begin{equation*}
    \Phi^{\epsilon}_i(x) = \frac{1}{2} \left( G_i(x) + H_i(x) -
    \sqrt{(G_i(x)-H_i(x))^2 + \epsilon^2} \right).
\end{equation*}
This function satisfies
\begin{equation*}
\begin{aligned}
    &\Phi^{\epsilon}_i(x) = 0 \quad \Longleftrightarrow \\
    &G_i(x) > 0, \; H_i(x) > 0, \; G_i(x)H_i(x) = \epsilon^2/4 .
\end{aligned}
\end{equation*}
This allows to reformulate the MPCC into the form
\begin{equation*}
\begin{aligned}
    {\rm NCP}(\epsilon): \; 
    \min \; &\varphi(x) &{\rm multipliers}\\
    {\rm s.t.} \; &g_I(x) \leq 0, &v^{g_I} \\
    ~ &g_{E}(x) = 0, &v^{g_E}\\
    ~ &\Phi^{\epsilon}_i(x) = 0, \; i = 1 \dots n_c. &v^{\Phi}_i
\end{aligned}
\end{equation*}
As $\epsilon \to 0$, the corresponding sequence of stationary points $x^k$ of NCP($\epsilon^k$), converges to a stationary point $\bar x$ of the MPCC. 
Assuming MPCC-LICQ at $\bar x$, convergence properties of this NCP-based scheme have been summarized in \cite{wang-biegleronline} based on the following Lagrange function
\begin{equation*}
\mathcal L = \varphi(x) 
+ (v^{g_I})^T g_I(x) 
+ (v^{g_E})^T g_E(x)
- (v^{\Phi})^T \Phi^{\epsilon}(x) .
\end{equation*}
Note that since the sign for the last term is arbitrary, it was chosen
to be consistent with the underlying constraints $G(x)>0, H(x)>0$.
The convergence properties have been shown to be similar to those of the REG method. In addition,
\citet{fukushima-pang} have established the convergence to B-stationarity for a NCP-based scheme, under MPCC-LICQ at an accumulation point $\bar x$ and the second-order necessary conditions at every stationary point
$x^k$, together with an \textit{asymptotic weak nondegeneracy condition}, i.e., for every $i \in I_{GH}(\bar x)$, $G_i(x^k)$ and $H_i(x^k)$ approach zero in the same order of magnitude.

So far the convergence results have shown that although convergence to B-stationarity is desirable, it is usually established under restrictive conditions.
In particular, MPCC-LICQ is assumed so that B-stationarity is equivalent to S-stationarity.
For the more general case, where a local minimum is B-stationary but not S-stationary, we are short of theory for NLP-based methods to converge to such a solution.  
On the other hand, numerical practice has seen the influence of such a local minimum on the performance of solution
strategies.
For example, REG($\epsilon$) may give rise to large NLP multipliers
$v^{REG}_i$ for $i \in I_{GH}(\bar x)$, and at the same time,
the convergence can be slow and inaccurate.
Related difficulties with this approach are described and analyzed in
\cite{laemmel-shikhman, wang-biegleronline}. 

Another complexity in convergence is that the solution strategies could converge to a spurious solution $\bar x$ at which $I_{GH}(\bar x) \neq \emptyset$.
Many MPCC problems having such a feasible but not optimal point have been presented in the literature; see, for instance, \cite{guo-lin-ye}, \cite{fletcher}, and \cite{kazi-thombre-biegler}. 

\subsection{Geometry Simplification}

In order to ensure first-order optimality, verifying the
B-stationarity condition (\ref{bstat}) directly is generally non-trivial, because the tangent cone $\mathcal T(\bar x)$ representing the geometry of the constraints at $\bar x$ is usually not easy to be expressed analytically. 

So we define B-stationarity in terms of the linearized tangent cone:
\begin{equation*} 
\begin{aligned}
    \mathcal T_{\rm MPCC}^{lin}&(\bar x) = \{ d| \\
    ~ & \nabla g_{I_i}(\bar x)^T d \leq 0, 
      &&\forall i \in I_{g_I}(\bar x),\\
    ~ & \nabla g_{E_i}(\bar x)^T d = 0, 
      &&\forall i = 1 \dots |g_E|,\\
    ~ & \nabla G_i(\bar x)^T d = 0, 
      &&\forall i \in I_{G}(\bar x),\\
    ~ & \nabla H_i(\bar x)^T d = 0, 
      &&\forall i \in I_{H}(\bar x),\\
    ~ & \nabla G_i(\bar x)^T d \geq 0, 
      &&\forall i \in I_{GH}(\bar x),\\
    ~ & \nabla H_i(\bar x)^T d \geq 0, 
      &&\forall i \in I_{GH}(\bar x),\\
    ~ & (\nabla G_i(\bar x)^T d) \cdot (\nabla H_i(\bar z)^T d) = 0,
      &&\forall i \in I_{GH}(\bar x) \} .
\end{aligned}
\end{equation*}
This is an analytic expression of a
simplified tangent cone and preserves the complementarity structure.
By assuming
\begin{equation}\label{mpcc-acq}
    \mathcal T(\bar x) = \mathcal T_{\rm MPCC}^{lin}(\bar x),
\end{equation}
the B-stationarity condition
(\ref{bstat}) is transformed into: 
\begin{equation} \label{bstat-acq}
    \nabla \varphi(\bar x)^T d \geq 0, \quad
    \forall d \in \mathcal T_{\rm MPCC}^{lin}(\bar x) ,
\end{equation}
meaning that at a B-stationary point $\bar x$, there are no descent
directions in the linearized tangent cone $\mathcal T_{\rm
  MPCC}^{lin}(\bar x)$.
Equivalently, $\bar x$ is a B-stationary point of the MPCC if and only if $d = 0$ solves the following linear program with complementarity constraints (LPCC):
\begin{equation}\label{lpcc}
\begin{aligned}
    \min \quad &\nabla \varphi(\bar x)^T d \\
    {\rm s.t.} \quad &\nabla g_{I}(\bar x)^T d \leq 0, \; \nabla g_E(\bar x)^T d = 0,\\
    ~ &\nabla G_{i}(\bar x)^T d = 0, \quad 
       \forall i \in I_{G}(\bar x), \\
    ~ &\nabla H_{i}(\bar x)^T d = 0, \quad
       \forall i \in I_{H}(\bar x), \\
    ~ &0 \leq \nabla G_{i}(\bar x)^T d \perp 
              \nabla H_{i}(\bar x)^T d \geq 0,
       \forall i \in I_{GH}(\bar x) .
\end{aligned}
\end{equation}
Note that (\ref{bstat-acq}) (equivalently, (\ref{lpcc})) is a stronger stationarity property than (\ref{bstat}) because it always holds that $\mathcal T_{\rm MPCC}^{lin}(\bar x) \supseteq \mathcal T(\bar x)$ (\cite{flegel-kanzowacq}).

\subsection{MPCC Algorithm}

This section develops an algorithm that explores the LPCC at a
feasible point of an MPCC and improves the objective whenever the
current point is not B-stationary.

For the MPCC problem, we consider solutions by the NCP($\epsilon$) or REG($\epsilon$) approach.
At the end of the loop, we arrive at a stationary point $\bar{x}$. To
check for B-stationarity of $\bar x$, we apply appropriate stationarity criteria to determine whether LPCC
(\ref{lpcc}) needs to be solved. In order to facilitate the solution of the LPCC,
we introduce binary variables to the complementarity constraints and
obtain the following mixed integer linear program (MILP):

\begin{equation} \label{eqn:milp}
\begin{aligned}
 \min_d \quad & \nabla \varphi(\bar{x})^T d \\
{\rm s.t.} \quad &\nabla g_{I_i}(\bar{x})^T d \leq 0, \quad i \in
I_{g_I}(\bar x), \\
&  \nabla g_E(\bar{x})^T d = 0, \\
& \nabla G_i(\bar{x})^T d = 0, \quad i \in I_G(\bar x), \\
& \nabla H_i(\bar{x})^T d = 0, \quad i \in I_H(\bar x), \\
& 0 \leq \nabla G_i(\bar{x})^T d \leq M w_i, \\
& 0 \leq \nabla H_i(\bar{x})^T d \leq M (1-w_i), \\
& w_i \in \{0, 1\},  \quad i \in I_{GH}(\bar x), 
\end{aligned}
\end{equation}
where $M$ is a positive constant. 
Note that B-stationary is satisfied if $d=0$ at the solution of (\ref{eqn:milp}).
Otherwise, we note that the MILP solutions $d/\|d\|$ and
binary variables $w_i$ remain unchanged for any positive value of $M$. 

If $d\neq 0$ solves (\ref{eqn:milp}), then we 
need to consider the following NLP with relaxed bi-active constraints based on binary variables from solution of the MILP (\ref{eqn:milp}):

\begin{equation}  \label{eqn:subnlp}
\begin{aligned}
\min_d \quad & \varphi(x) \\  
{\rm s.t.} \quad & \varphi(x) \leq \varphi(\bar{x}), \\ 
& g_I(x) \leq 0, \; g_E(x) = 0,\\
& G_i(x) = 0, \; H_i(x) \geq 0, \quad i \in I_G(\bar x), \\
& G_i(x) \geq 0, \; H_i(x) = 0, \quad i \in I_H(\bar x), \\
& G_i(x) \geq 0, \; H_i(x) \geq 0, \quad i \in I_{GH}(\bar x), \\
& G_i(x) = 0, \; H_i(x) \geq 0, \quad \mbox{if $w_i = 0$}, \\
& G_i(x) \geq 0, \; H_i(x) = 0, \quad \mbox{if $w_i = 1$} .
\end{aligned}
\end{equation}

Using these problem formulations, we now state our MPCC Algorithm. 

\renewcommand{\thealgorithm}{} 
\begin{algorithm}\label{mpcc-b-algorithm}
\caption{MPCC algorithm for B-stationarity}
\begin{algorithmic}[0] 
\setlength{\itemsep}{5pt} 
\STATE \textbf{Step 1:} 
Solve NCP($\epsilon$) or REG($\epsilon$) with $\epsilon^k \to 0$ (We caution that  REG($\epsilon$) may struggle in converging to solutions that are not S-stationary.) 
\STATE \textbf{Step 2:}
At the solution $\bar{x}$ of NCP($\epsilon$) or REG($\epsilon$):
\begin{itemize}
\item If $I_{GH}(\bar x) = \emptyset$, STOP. 
There are no bi-active constraints and  $\bar{x}$ is B-stationary. 
\item For NCP($\epsilon$), if $v^{\Phi}_i  \geq  0$ for all $i \in I_{GH}(\bar x)$, STOP. For REG($\epsilon$), if $v_i^{REG}$ are well-scaled and appear bounded for all $i \in I_{GH}(\bar x)$, STOP.
These conditions are equivalent to $\sigma^G_i, \sigma^H_i \geq 0$ for all $i \in I_{GH}(\bar x)$ and $\bar{x}$ is S-stationary (see \cite{wang-biegleronline, scholtes2001} for details).
\end{itemize} 
\STATE \textbf{Step 3:}
At $\bar{x}$, solve MILP (\ref{eqn:milp}) to attempt to find a descent direction $d$ from the bi-active complementarities,
represented by binary variables $w$. If $d = 0$, STOP and $\bar{x}$ is B-stationary. 
\STATE \textbf{Step 4:} 
Starting from $\bar{x}$, solve NLP (\ref{eqn:subnlp}). Designate this NLP solution 
as $\bar{x}$. If $I_{GH}(\bar x) = \emptyset$, STOP as $\bar x$ is B-stationary; otherwise, go to Step 3. 
\end{algorithmic}
\end{algorithm}

%
%
%

This MPCC algorithm is patterned after the MPCC solution strategy proposed in \cite{nurkanovic2025hybrid}. In addition, our approach includes the following modifications and simplifications as well, which lead to smaller LPCCs and leverage the performance of large-scale NLP solvers. 
\begin{itemize} 
\item The NLPs, NCP($\epsilon$), REG($\epsilon$) and (\ref{eqn:subnlp}), which represent reformulations of the MPCC, are fully converged to KKT points. These are feasible points for the MPCC and allow identification of the sets
$I_{G}, I_{H}$, and $I_{GH}$.  The NCP($\epsilon$) approach may be preferred as this approach allows convergence to C-, M-, as well as S- stationary points (\cite{wang-biegleronline}).

\item As in \cite{nurkanovic2025hybrid}, an LPCC is formulated into an MILP to determine whether the NLP
solutions are B-stationary. However, the MILP (\ref{eqn:milp}) uses binary variables 
only for the bi-active complementarity constraints, and often leads to a significant reduction 
in binary variables for large MPCCs, e.g., derived from discretized dynamic models. 

\item The MILP (\ref{eqn:milp}) is written in a simple form that determines relaxed bi-active constraints corresponding to descent directions. The optimal relaxations, determined by $w$, are independent of $M>0$
and eliminate the need for a trust region to determine potential descent directions. 

\item As in \cite{nurkanovic2025hybrid}, only the binary variables $w$
  are used from the solution of (\ref{eqn:milp}), in order to formulate the constraints in (\ref{eqn:subnlp}) and determine an improved point for the MPCC. This new point is then evaluated for B-stationarity in the next cycle. 


\item The algorithm is scalable to large problems, in the sense that our subproblems are solved with NLPs that efficiently handle large-scale systems. In addition, to check a given point for B-stationarity that is not S-stationarity, we only need to consider a small 
MILP. Moreover, neither determining descent directions (Step 3) nor finding a new trial point with improved objective (Step 4) requires expensive trust-region techniques to confine the feasible region. Finally, for hybrid dynamical systems, we consider only the switching points identified by bi-active complementarities to determine the size of the MILP, instead of including all potential switching points. This greatly reduces the MILP problem size, and therefore requires far less computational effort.
\end{itemize}

\citet{nurkanovic2025hybrid} have proved that their solution strategy converges to a point that is either locally infeasible or B-stationary for MPCC \eqref{eq:general-MPCC}, under the following assumptions: 
\begin{itemize}
\item An NLP solver applied to the NLPs (REG($\epsilon$) and (\ref{eqn:subnlp})) with starting point $x_0$, either returns a stationary point $\bar{x}$ or a certificate of local infeasibility.

\item The above MILP problem applied at feasible $\bar{x}$ returns either a descent direction with $\nabla \varphi(x)^Td < 0$ or $d=0$.

\item There exists a compact set that contains the feasible set of the MPCC.

\item The condition \eqref{mpcc-acq} holds at every feasible point of the MPCC. 
\end{itemize}
Because our algorithmic simplifications retain the essential elements of their approach, their convergence proofs also apply to the algorithm developed in this section.

\subsection{Illustrative Examples for MPCC Algorithm} \label{illustrate}

To illustrate the proposed algorithm, we present a few examples.

\subsubsection{Example 1 (M-stationary local minimum).}
Consider the following problem from \cite{guo-lin-ye}:
\begin{equation*}
\begin{aligned}
  \min \quad &x_1 + x_2 - x_3 \; &{\rm multipliers}\\
  {\rm s.t.} \quad
     &-4x_1 + x_3 \leq 0, &\lambda_1\\
  ~  &-4x_2 + x_3 \leq 0, &\lambda_2\\
  ~  &0 \leq x_1 \perp x_2 \geq 0 . &\sigma_1,\sigma_2
\end{aligned}
\end{equation*}
The unique minimum $\bar x =(0, 0, 0)$ is not S-stationary, due to the multipliers $\sigma_1+\sigma_2 = -2$. 

Step 1: Apply IPOPT or CONOPT to solve a sequence of problems NCP($\epsilon$). 

Step 2: Consider the solution of NCP($\epsilon$) at $\bar{x} = (0, 0, 0)$. With the weak stationarity conditions given by:
\begin{equation*}
\begin{aligned}
1 - 4 \lambda_1 - \sigma_1 &= 0, \\
1 - 4 \lambda_2 - \sigma_2 &= 0, \\
\lambda_1 + \lambda_2 &= 1, \\
\end{aligned}
\end{equation*}
the multipliers $ \sigma_1 + \sigma_2 =-2$, and the solution is at best M-stationary.

Step 3: Formulate the corresponding MILP (\ref{eqn:milp}):
\begin{equation*}
\begin{aligned}
\min \quad & d_1 + d_2 - d_3 \\
{\rm s.t.} \quad & -4d_1 + d_3 \leq 0, \\
& -4d_2 + d_3 \leq 0, \\
& 0 \leq d_1 \leq M w, \\
& 0 \leq d_2 \leq M (1-w) , \quad w \in \{0, 1\},
\end{aligned}
\end{equation*}
which has the solution $(d_1, d_2, d_3) = (0, 0, 0)$, 
confirming that $\bar{x}$ is a B-stationary point. 

\subsubsection{Example 2 (Nonoptimal M-stationary point).}
Consider the following problem:
\begin{equation*}
\begin{aligned}
  \min \quad &(x_1 - 1)^2 + x_2^2 \; &{\rm multipliers}\\
  {\rm s.t.} \quad
     &x_1 \leq 1, &\lambda_1\\
  ~  &x_2 \geq 0, &\lambda_2\\
  ~  &0 \leq x_1 \perp x_2 \geq 0 . &\sigma_1,\sigma_2
\end{aligned}
\end{equation*}
The global minimum is $(1,0)$, which is S-stationary.
Consider the M-stationary point $\bar x =(0,0)$, where the multipliers are $(\sigma_1,\sigma_2) = (-2,0)$.
We find that $\bar x$ is not B-stationary because $d = 0$ is not optimal to the LPCC (\ref{lpcc}); every feasible direction $d =(d_1>0, d_2=0)$ leads to a reduction of the objective $\nabla \varphi(\bar x)^T d = -2d_1 <0$. 

Step 1: Apply IPOPT or CONOPT to solve a sequence of problems NCP($\epsilon$). When initialized close to the M-stationary point $\bar x = (0,0)$, both solvers move away and converge to the global minimum.
Unfortunately, when initialized exactly at the M-stationary point, CONOPT declares a local optimum at this point; this is likely due to lack of precision at $(0,0)$ in finding an appropriate active set. 

Step 2: Consider the solution of NCP($\epsilon$) at $\bar{x} = (0, 0)$. The weak stationarity conditions lead to $(\sigma_1, \sigma_2) = (-2, 0)$, so that $\bar{x}$ is M-stationary.

Step 3: Formulate the corresponding MILP (\ref{eqn:milp}):
\begin{equation*}
\begin{aligned}
\min \quad & -2 d_1  \\
{\rm s.t.} \quad & 0 \leq d_1 \leq M w, \\
& 0 \leq d_2 \leq M (1-w) , \quad w \in \{0, 1\},
\end{aligned}
\end{equation*}
which has the solution $(d_1, d_2) = (M, 0)$ and $w = 1$, providing a descent direction.

Step 4: Solve the relaxed NLP (\ref{eqn:subnlp}): 
\begin{equation*}
\begin{aligned}
  \min \quad &(x_1 - 1)^2 + x_2^2 \\
  {\rm s.t.} \quad &(x_1 - 1)^2 + x_2^2 \leq 1, \\
  &x_1 \leq 1, \\
  &x_2 \geq 0, \\
  &x_1 \geq 0, \; x_2 = 0,
\end{aligned}
\end{equation*}
which yields the solution $\bar{x} = (1, 0)$. The MPCC has no bi-active complementary constraints at $\bar x$ and thus it is a B-stationary point.
  
\subsubsection{Example 3 (C-stationary local maximum).}
Consider the following problem from \cite{kazi-thombre-biegler}:
\begin{equation*}
\begin{aligned}
  \min \quad &(x_1-1)^2 + (x_2-1)^2 \; &{\rm multipliers}\\
  {\rm s.t.} \quad
     &0 \leq x_1 \perp x_2 \geq 0 . &\sigma_1,\sigma_2
\end{aligned}
\end{equation*}
This problem has two B-stationary points, $(1,0)$ and $(0,1)$, which are global minima. There is also a local maximum and C-stationary point at $\bar x = (0,0)$, where the multipliers are $(\sigma_1,\sigma_2) = (-2,-2)$.
The point $\bar x$ is not B-stationary because $d = 0$ does not solve the LPCC (\ref{lpcc}); every feasible direction $d$ with $d_1 >0$ or $d_2 >0$ leads to a reduction $\nabla \varphi(\bar x)^T d <0$.

Step 1: For the NCP formulation starting from the C-stationary point $\bar x = (0,0)$, CONOPT converges to one of the minima, depending on the value of $\epsilon$.  On the other hand, IPOPT remains at $(0,0)$ for small initial values of $\epsilon$ because IPOPT strongly regularizes the reduced Hessian of the Lagrange function, which is negative definite. Similarly, SNOPT also remains at
$(0,0)$ for small initial values of $\epsilon$, because it does not use second derivatives but a BFGS approximation to the Hessian instead.

Step 2: Consider the solution of NCP($\epsilon$) at $\bar{x} = (0, 0)$. The weak stationarity conditions lead to $(\sigma_1, \sigma_2) = (-2, -2)$ and $\bar{x}$ is C-stationary.

Step 3: Formulate the corresponding MILP (\ref{eqn:milp}):
\begin{equation*}
\begin{aligned}
\min \quad & -2 d_1 - 2 d_1  \\
{\rm s.t.} \quad & 0 \leq d_1 \leq M w, \\
& 0 \leq d_2 \leq M (1-w) , \quad w \in \{0, 1\},
\end{aligned}
\end{equation*}
which has the solution $(d_1, d_2) = (M, 0)$ and $w = 1$, or $(d_1, d_2) = (0,M)$ and $w = 0$, either of which provides a descent direction.

Step 4: Solve the relaxed NLP (\ref{eqn:subnlp}): 
\begin{equation*}
\begin{aligned}
  \min \quad &(x_1 - 1)^2 + (x_2-1)^2 \\
  {\rm s.t.} \quad  &(x_1 - 1)^2 + (x_2-1)^2 \leq 2,\\
  &x_1 \geq 0, \; x_2 = 0, \quad   
  ( \mbox{or } x_1 = 0, \; x_2 \geq 0 ) 
\end{aligned}
\end{equation*}
which yields the solution $\bar{x} = (1, 0)$ (or $\bar{x} = (0, 1)$). The MPCC has no bi-active constraints at $\bar x$ and thus it is a B-stationary point.

\section{DCS Numerical Studies}

In this section, we solve three examples from literature to showcase the implementation and efficacy of our proposed strategy. The first example is a simple hybrid dynamical system with \textit{signum} function. The system has a closed form piecewise solution with exactly one switch point. Our method accurately identifies the switch point and the non-unique, non-uniform discretization. The second example considers an optimal control problem with multiple variables and switching functions. Using our proposed strategy, we obtain the optimal sliding mode solution and accurate location of multiple switching points. Lastly, the third case study is based on an engineering example with the optimal control of a gas-liquid tank system. This example involves controlling the valve opening to minimize the final liquid level and is inspired by \citet{moudgalya2001class-I}. In this example, we demonstrate that our method is able to solve hybrid dynamical problems with a large number of finite elements efficiently and able to accurately locate the switch points.

For implementation, we use Julia's optimization package JuMP to formulate the optimization problems and solve it using IPOPT as the default NLP solver. In every example, we first discretize the hybrid dynamical system using first-order implicit Euler method and the moving finite element formulation Eq.\eqref{eq:cross-complementarity}. Next we reformulate the MPCC using either the REG($\epsilon$) or NCP($\epsilon$) and solve the NLP with the relaxation parameter $\epsilon=0.1$. We then re-discretize the differential equations using either Runge-Kutta (RK) scheme or Radau based orthogonal collocation method (OCM) for higher order accuracy. We initialize the variables using the implicit Euler solution. We then solve this large discretized NLP recursively by initializing it with the previous NLP solution and with smaller values of $\epsilon$ till it reaches $\epsilon_{tol} = 10^{-6}$. Then, we calculate the location of switch points using the Eq.\eqref{eq:switch-point} and add the step-equilibration constraints Eq.\eqref{eq:uniform_step} and Eq.\eqref{eq:non-uniform_step} to the discretized NLP. Finally, we solve the NLP with higher order scheme, cross-complementarity constraints (\ref{eq:cross-complementarity}f)
and step-equilibration constraints \eqref{eq:proposed-formulation} to find the optimal solution. We then determine the bi-active set ($I_{GH}$) of the solution and implement the MPCC algorithm described in Section \ref{mpcc-b-algorithm} to check whether we obtain a B-stationary solution to the optimization problem.

All problems were solved to a tolerance = $10^{-8}$ within 20 CPU seconds using a standard M3-MacBook Air machine. Non-default IPOPT options with ``mu\_strategy - adaptive" and ``linear\_solver - ma27" were used to improve the convergence of the NLPs. 

\subsection{Signum Problem}

The signum problem is a commonly used example for hybrid dynamic systems (presented below). 
\begin{subequations}\label{eq:sign}
    \begin{align}
        \min_{x} \enskip & (x(t_f) - 5/3)^2, \qquad  t \in [0,2],\\
        \text{s.t.} \enskip & \dot x = 2 - \text{sgn}(x), \quad x(0) = -2.
    \end{align}
    Although the problem has zero degrees of freedom, the problem is formulated as an optimization problem for illustrative purposes. We implement three different strategies to solve the sign problem: (1) Moving Finite Element (MFE) strategy proposed in \citet{Baumrucker2009mpec} (2) Step Equilibration strategy in \citet{nurkanovic2022finite} and (3) The proposed two stage step equilibration strategy introduced in this study. 
    
    The problem is solved by discretizing the time interval into $N=10$ finite elements. We implement the 4th order RK4 discretization scheme to transform the differential equation to algebraic constraints. 
    
    The sgn($x$) is described as a piecewise function:
    \begin{equation}
       \text{sgn}(x)  =  \begin{cases} 
      -1 & x < 0, \\
      [-1,1] & x = 0,\\
      1 &  x > 0 .
   \end{cases}
    \end{equation}
    The piecewise formulation is reformulated as a DCS \eqref{eq:DCS} using the complementarity constraints as:
    \begin{align}
        & \dot x = 3 \nu_1 + \nu_2,\\
        & x - \lambda_1 - \mu = 0, \\
        & -x - \lambda_2 - \mu = 0, \\
        & \nu_1 + \nu_2 = 1,\\
        & 0 \leq \nu_i \perp \lambda_i \geq 0, \enskip i = 1,2 .
    \end{align}
    We can substitute $\nu_2 = 1 - \nu_1$ and subtract (\ref{eq:sign}e) and (\ref{eq:sign}f) to reduce the system as:
    \begin{align}
        & \dot x = 1 + 2\nu,\\
        & 2x - \lambda_1 + \lambda_2 = 0,\\
        & 0 \leq \nu \perp \lambda_1 \geq 0,\\
        & 0 \leq 1-\nu \perp \lambda_2 \geq 0 .
    \end{align}
\end{subequations}

The analytical solution to the problem in Eq.\ref{eq:sign} is given by:  
\begin{equation}\label{eq:accurate}
    x^* = 
    \begin{cases}
    -2 + 3t, & \text{if } t < 2/3, \\
    t-2/3, & \text{if } t \geq 2/3.
    \end{cases}
\end{equation}

\begin{figure}[H]
    \centering
    \includegraphics[width=0.75\linewidth]{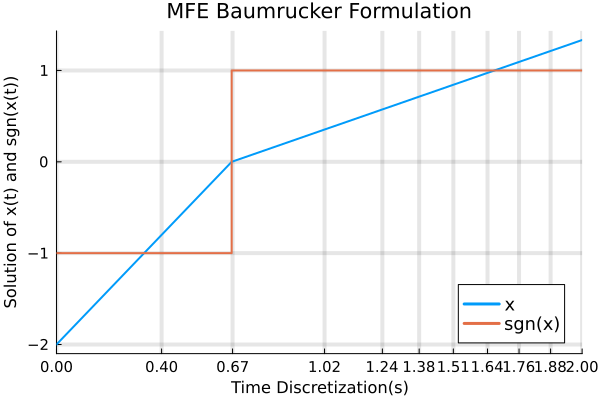}
    \caption{Solution profile for state $x$ and sign variable using the MFE strategy in \cite{Baumrucker2009mpec}}
    \label{fig:signum_plot_brian}
\end{figure}

\begin{figure}[H]
    \centering
    \includegraphics[width=0.75\linewidth]{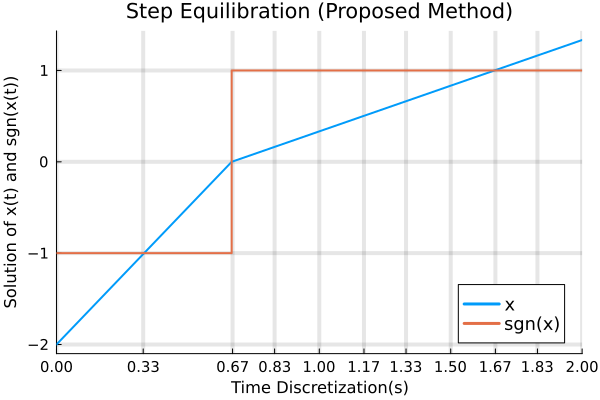}
    \caption{Solution profile for state $x$ and sign variable using the proposed two step MFE strategy}
    \label{fig:signum_plot_proposed}
\end{figure}

The numerical solution for $x$ and the sign variable $1-2\nu$ is plotted in Fig.\ref{fig:signum_plot_brian} and Fig.\ref{fig:signum_plot_proposed} for MFE strategy by \cite{Baumrucker2009mpec} (Baumrucker et al.) and the proposed two stage step equilibration strategy respectively. Both the state variable $x$ and switch variable $v$ is accurate approximation of the analytical solution Eq.\eqref{eq:accurate} for all methods. As can be seen in the Fig.\ref{fig:signum_plot_brian} and Fig.\ref{fig:signum_plot_proposed}, the switch point and switching time $t_s = 2/3$ is accurately located in both cases and ensured that the switch occurs at the element boundary. 

But the two methods result in different discretization and step size values. In Fig.\ref{fig:signum_plot_brian}, the Baumrucker's MFE strategy produces non-uniform and non-unique discretization results $h_1 = 0.40, h_2 = 0.27, h_3 = 0.35, h_4 = 0.22, h_5 = 0.14, h_6 = h_7 = 0.13$ and $h_8 = h_9 = h_{10} = 0.12$. In comparison, our proposed two-stage approach produces a uniform and unique discretization, except at the switch point $h_1 = h_2 = 0.333, h_3 = h_4 = h_5 = h_6 = h_7 = h_8 = h_9 = h_{10} = 0.167$.   
\vspace{-3mm}
Next, we compare our proposed method with the Nurkanovic formulation for different number of finite elements $N=10,20,50$ and $100$. While both methods are able to solve all the instances to same solution, it is evident from Table \ref{tab:ipopt_iterations} that it takes more computational effort (\# of IPOPT iterations) to solve the problem with the Nurkanovic formulation as compared to the two step approach. 
\begin{table}[h]
    \centering
    \begin{tabular}{|c|c|c|}
    \hline
    No. of elements & Proposed Method & Nurkanovic\\
    \hline
    10 & 26 & 470 \\
    20 & 56 & 759 \\
    50 & 65 & 750 \\
    100 & 303 & 1162\\
    \hline
    \end{tabular}
    \caption{Number of iterations in IPOPT}
    \label{tab:ipopt_iterations}
\end{table}
\vspace{-3mm}
We also observe that the solution obtained using the Nurkanovic formulation exhibited some numerical instability for values of variable $\tau$ ranging from $10^{-9}$ to $1000-3500$. This results in a lot of iterations going to the restoration phase in IPOPT, suggesting that the formulation and the resulting problem is ill-conditioned.

Note that although the Nurkanovic formulation implemented and described in this study in section 2.3 may have a few differences in implementation from the NOSNOC solver \cite{nurkanovic2022nosnoc}, we expect the general trend to be similar. Finally, the optimal solution for this example has an empty bi-active set and thus is a B-stationary point.


\subsection{Sign Optimal Control Problem}

The second example, taken from \cite{nurkanovic2022finite}, involves an optimal control problem with 4 state variables $(x_1,x_2,x_3,x_4)$ and 2 control inputs $(u_1,u_2)$. The aim of the OCP is to drive the states $(x_1,x_2)$ to the target point $(x^s_1,x^s_2) = (-\pi/6,-\pi/4)$ with minimum control input and secondary inputs $(x_3,x_4)$ as seen in the objective function.
\begin{subequations}\label{eq:sign_ocp}
    \begin{align}
        \min_{u} \enskip & \int^4_0 (u(t)^T u(t) + x^2_3(t) + x^2_4(t))dt \nonumber \\
        & + \rho ((x_1(4) - x^s_1)^2 +  (x_2(4) - x^s_2)^2)  \\
        \text{s.t.} \enskip & \quad  \mathbf{x}(0) = (2\pi/3,\pi/3,0,0), \\
        & \quad \dot x_1(t) = -\text{sgn}(\psi_1(x)) + x_3(t), \\
        & \quad \dot x_2(t) = -\text{sgn}(\psi_2(x)) + x_4(t), \\
        & \quad \dot x_3(t) = u_1(t), \\
        & \quad \dot x_4(t) = u_2(t), \\
        &  -2 \leq x_3(t), x_4(t) \leq 2, \\
        &  -10 \leq u_1(t), u_2(t) \leq 10,
    \end{align}
where $\psi_1(x) = x_1 + 0.15x_2^2$ and $\psi_2(x) = -0.05x_1^3 + x_2$ and $\rho = 10^3.$ 
\end{subequations}
There are two switching functions $\psi_1(x)$ and $\psi_2(x)$. The sign-based switch functions are reformulated using auxiliary slack variables $(s^p_1,s^n_1,s^p_2,s^n_2)$ and multiplier variables $(\nu_1,\nu_2)$ as shown in Eq.\eqref{eq:sign_ocp_ref}:
\begin{subequations}\label{eq:sign_ocp_ref}
    \begin{align}
        & \dot x_1(t) = - \nu_1(t) + x_3(t), \\
        & \dot x_2(t) = - \nu_2(t) + x_4(t), \\
        & \psi_1(x) = s^p_1(t) - s^n_1(t), \\
        & \psi_2(x) = s^p_2(t) - s^n_2(t), \\
        & 0 \leq 1-\nu_1(t) \perp s^p_1(t) \geq 0, \\
        & 0 \leq 1+\nu_1(t) \perp s^n_1(t) \geq 0, \\
        & 0 \leq 1-\nu_2(t) \perp s^p_2(t) \geq 0, \\
        & 0 \leq 1+\nu_2(t) \perp s^n_2(t) \geq 0.
    \end{align}
\end{subequations}

\begin{figure}[H]
    \centering
    \includegraphics[width=0.8\linewidth]{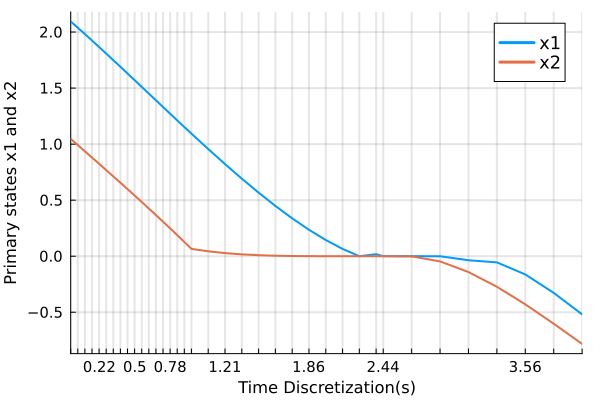}
    \caption{Plot of $x_1(t)$ and $x_2(t)$ vs $t$}
    \label{fig:sgn_ocp_x1_x2}
\end{figure}

\begin{figure}[H]
    \centering
    \includegraphics[width=0.8\linewidth]{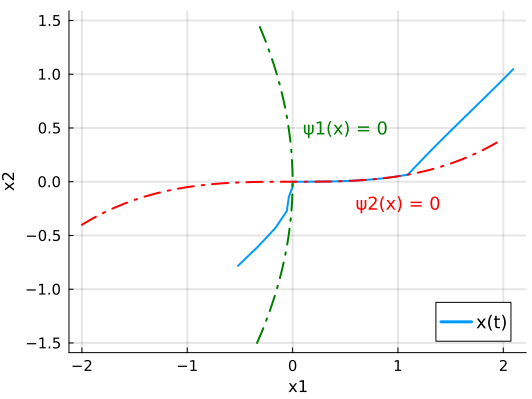}
    \caption{Plot of $x_1(t)$ vs $x_2(t)$ with $\psi_1 = \psi_2 = 0$ contour plots}
    \label{fig:sgn_ocp_x1_v_x2}
\end{figure}

The differential equations are discretized with $N=36$ finite elements for $t_f = 4$ using 3pt. Radau orthogonal collocation method (OCM). Figures \ref{fig:sgn_ocp_x1_x2} and \ref{fig:sgn_ocp_x1_v_x2} represent the optimal values of the state variables ($x_1,x_2$). As seen in Fig. \ref{fig:sgn_ocp_x1_x2}, there are four switch points at 17th, 28th, 29th and 32nd finite elements corresponding to $t = 0.944, 2.38, 2.44$ and $3.11$s respectively.  

Fig. \ref{fig:sgn_ocp_x1_v_x2} shows sliding modes on non-linear manifolds. The solution trajectory moves from: (1) $x(0)$ to $\psi_2(x) = 0$;
(2) slides on the curve to origin $\psi_2(x) = \psi_1(x) = 0$; (3)
stay at the origin and (4) slide on the curve $\psi_1(x) = 0$
before leaving it to reach $(x^s_1,x^s_2)$. This simple example
shows that our method can handle multiple switching
points and sliding modes on non-linear manifolds.

\begin{figure}[H]
    \centering
    \includegraphics[width=0.8\linewidth]{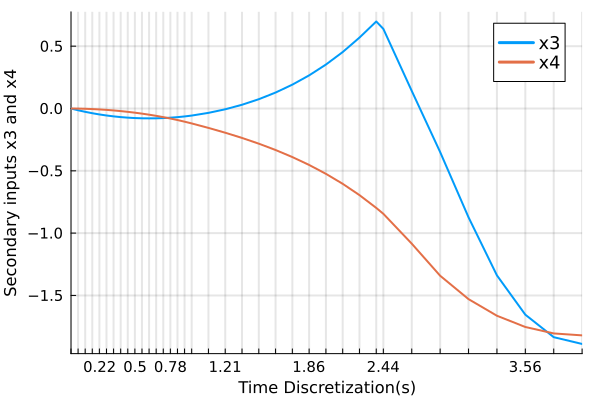}
    \caption{Plot of $x_3(t)$ and $x_4(t)$ vs $t$}
    \label{fig:sgn_ocp_x3_x4}
\end{figure}

\begin{figure}[H]
    \centering
    \includegraphics[width=0.8\linewidth]{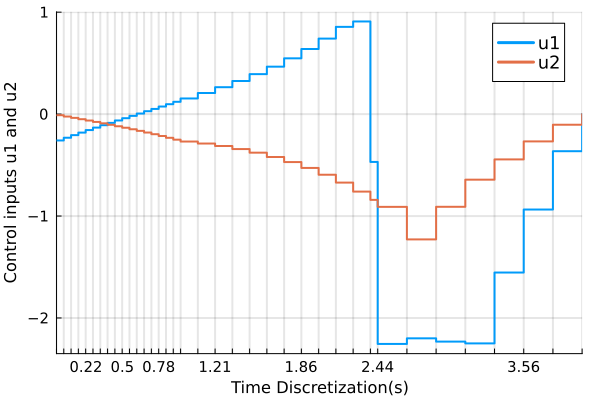}
    \caption{Plot of $u_1(t)$ and $u_2(t)$ controls vs $t$}
    \label{fig:sgn_ocp_u1_u2}
\end{figure}

 Figures \ref{fig:sgn_ocp_x3_x4} and \ref{fig:sgn_ocp_u1_u2} plot the profile of secondary state variables ($x_3,x_4$) and control inputs ($u_1,u_2$) respectively. Finally, the optimal solution has a non-empty bi-active set $I_{GH}$ at the 28th and 32nd finite element. Checking the bi-active multipliers $v^{REG}$ at the two finite elements (Step 2 of Algorithm), show that they are both appropriately scaled and bounded with $v^{max} = 17.72$ , which ensures that the solution is S-stationary and thus B-stationary.  

\subsection{Gas-Liquid Tank Problem}

\begin{figure}[h]
    \centering
    \includegraphics[scale = 0.3]{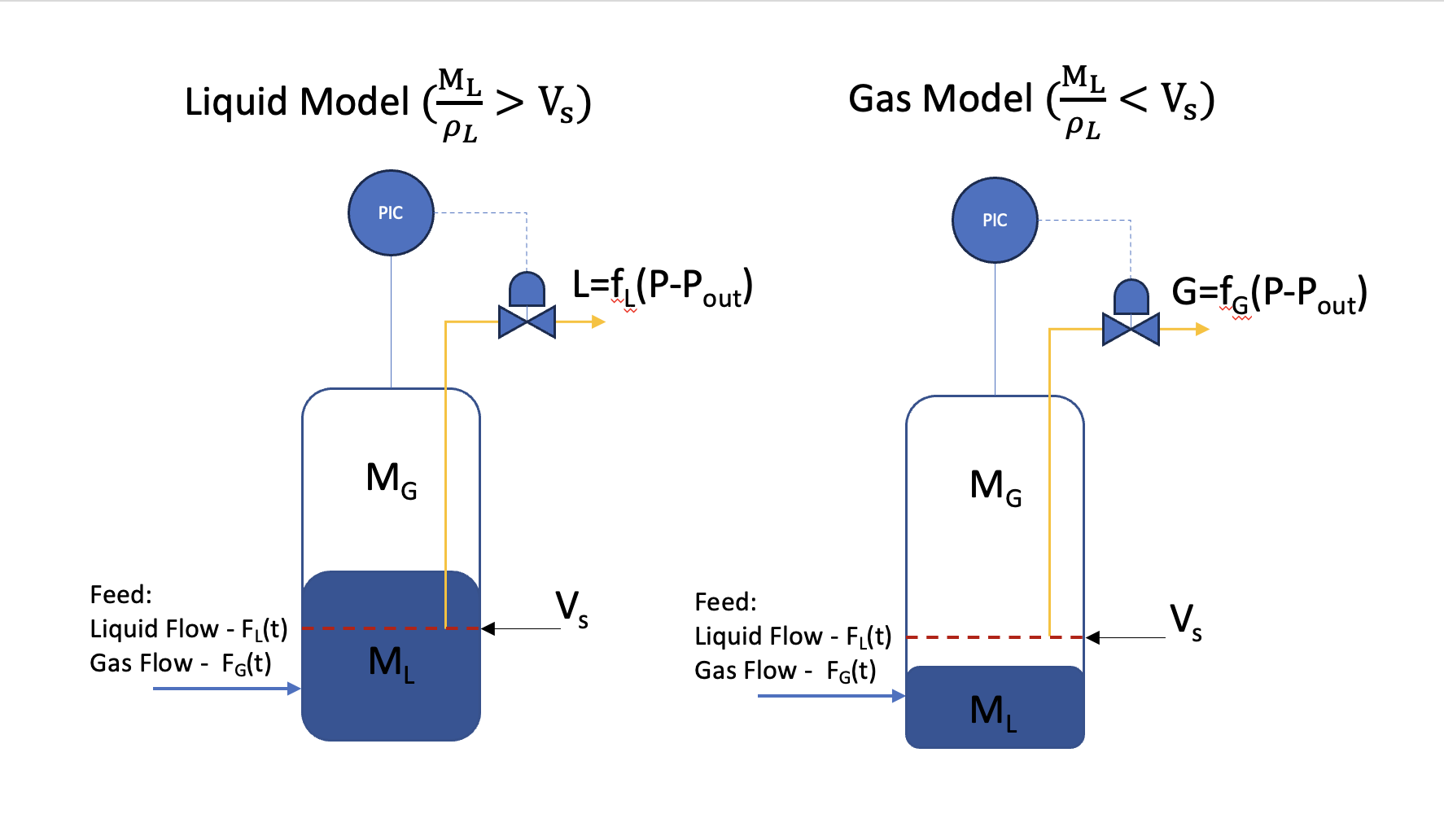}
    \caption{An ideal gas-liquid tank system with a pressure control valve}
    \label{fig:ideal gas-liquid system}
\end{figure}

The final case study with ideal gas-liquid tank system is an engineering example from \cite{moudgalya2001class-I}. In this example, there is a closed tank with one feed inlet and one outlet with a control valve that regulates the pressure inside the vessel as shown in Fig. \ref{fig:ideal gas-liquid system}.

The feed consists of a mixture of liquid ($F_L$) and ideal gas ($F_G$) and the outlet is either liquid ($L$) or gas ($G$) depending if the liquid level is above or below the outlet tube opening ($V_s$). The assumptions in the model are:
\begin{itemize}
    \item The gas and the liquid do not react.
    \item The liquid has negligible vapor pressure.
    \item The valve dynamics are ignored.
    \item The flow rate through the valve is proportional to the difference of tank and outlet pressure.
    \item The temperature, feed flow rates, outlet pressure and the valve opening are kept constant.
\end{itemize}

The dynamics of the system are described by the following set of differential equations.

\begin{subequations}\label{eq:liquid-gas}
    Liquid model: $(\frac{M_L}{\rho_L} > V_s)$ \hspace{5mm} Gas model: $(\frac{M_L}{\rho_L} < V_s)$
    \begin{align}
       & \frac{dM_G}{dt} = F_G, \hspace{20mm} && \frac{dM_G}{dt} = F_G - G, \\
       & \frac{dM_L}{dt} = F_L - L, \hspace{20mm} && \frac{dM_L}{dt} = F_L,\\
       & \qquad \qquad M_G\frac{RT}{P} + \frac{M_L}{\rho_L} = V , \\
       & L = k_L x (P - P_{out}). \hspace{10mm} && G = k_G x (P - P_{out}).
    \end{align}
\end{subequations}
The rates of liquid ($M_L$) and vapor holdup ($M_G$) vary according to the conservation based differential equations (\ref{eq:liquid-gas}a) and (\ref{eq:liquid-gas}b).  The total volume ($V$) and the tank pressure ($P$) are related by the ideal gas equation and liquid volume (\ref{eq:liquid-gas}c). The liquid and vapor outlet flow are related as a function of difference between tank and outlet pressure. The system parameter values in this example are listed in Table \ref{tab:params}.

We reformulate the above dynamic system into a DCS using complementarity constraints as:
\begin{subequations}
    \begin{align}
        & \frac{dM_G}{dt} = F_G \nu + (F_G - G)(1 - \nu), \\
        & \frac{dM_L}{dt} = (F_L - L) \nu + F_L(1 - \nu), \\
        & \frac{M_L}{\rho_L} - V_s = s_1 - s_2,\\
        & 0 \leq 1 - \nu \perp s_1 \geq 0,\\
        & 0 \leq \nu \perp s_2 \geq 0.
    \end{align}
    Here $s_1$ and $s_2$ are positive slack variables which denote the positive and negative part of the switching variable ($\frac{M_L}{\rho_L} - V_s$) respectively.
\end{subequations}
The objective function for this problem is a function of final liquid holdup and the control valve variable $x$ with the cost parameter $\beta = 100$
\begin{equation}
    \min \enskip M_L (t_f) + \beta \int^{25.0}_0 (x - 0.1)^2 dt .     
\end{equation}
\begin{table}[h!]
    \centering
    \begin{tabular}{|c|c|}
    \hline
    Parameters & Value \\
    \hline
    $F_L, F_G$(mol/sec) & 2.5, 0.1  \\
    $V, V_s$(litres) & 10, 5\\
    $T$(K)  & 300 \\
    $P_{out}$(atm) & 1 \\
    $\rho_L$(mol/l) & 50 \\
    $k_L, k_G$ & 1.0 \\
    \hline
    \end{tabular}
    \vspace{2mm}
    \caption{Parameter values in ideal gas-liquid closed system}
    \label{tab:params}
\end{table}
The initial conditions for the state variables are specified as: $P_0 = 35$ atm, $M_{L_0} = 260$ mol, $M_{G_0} = 6.83$ mol, $L_0 = G_0 = 0.25$ mol. Similar to the first example, the dynamic constraints are discretized with $N=100$ finite elements for a time horizon $t_f = 25$s. 

\begin{figure}[H]
    \centering
    \includegraphics[width=0.8\linewidth]{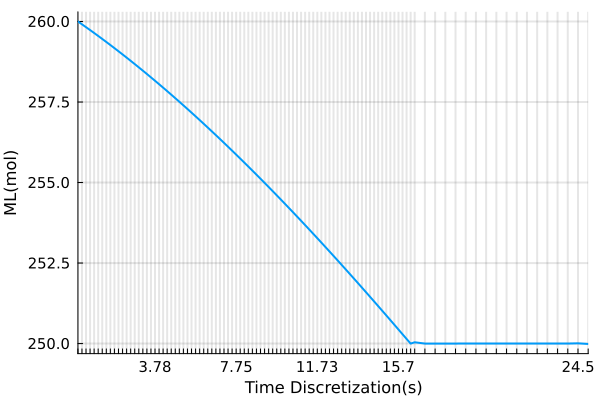}
    \caption{Profile of liquid holdup in the tank}
    \label{fig:liquid holdup}
\end{figure}

\begin{figure}[H]
    \centering
    \includegraphics[width=0.8\linewidth]{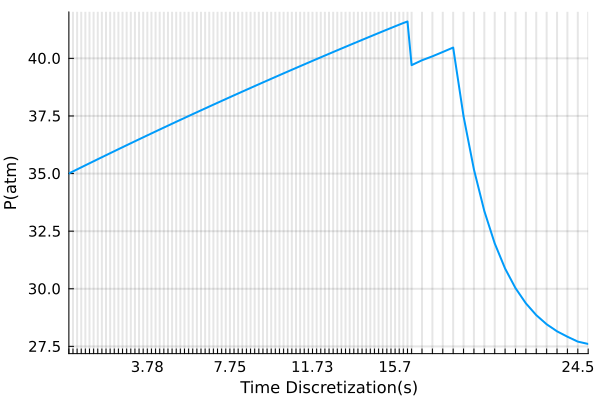}
    \caption{Pressure inside the tank vs time}
    \label{fig:tank pressure}
\end{figure}

The results from the optimization are plotted in Figures \ref{fig:liquid holdup} and \ref{fig:tank pressure}. The switching time or points in this case ($t_s = 16.3,16.5,18.5$s) can be clearly seen in the pressure plot. Initially, the system only has liquid exit as the liquid holdup is higher than the set point ($V_s$), the switch happens when the liquid holdup ($M_L = 250, M_L/\rho_L = 5.0$) goes below the set volume ($V_s$) and gas starts coming out from the outlet reducing the tank pressure as shown in the pressure profile Fig.\ref{fig:tank pressure}. 

\begin{figure}[H]
    \centering
    \includegraphics[width=0.8\linewidth]{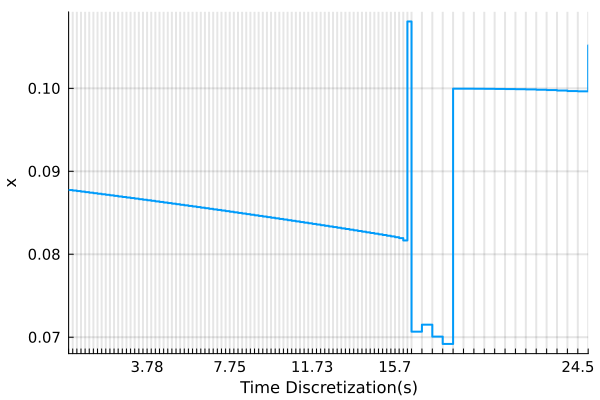}
    \caption{Control valve profile vs time}
    \label{fig:valve control variable}
\end{figure}

The plots clearly show that our method is able to identify the switching point in the hybrid dynamics and adapt the time step sizes ($h_l$) such that the switch time $t_s$ occurs at the boundary of an element.
As seen in Fig.\ref{fig:valve control variable}, the control variable $x$ varies almost linearly from $t=0$ to $t=16.3$s before abruptly increasing, then suddenly decreasing at $t=16.5$s and then increasing again at $t=18.5$s to a steady value of 0.1. The optimal solution has a non-zero bi-active set $I_{GH}$ active at the 83rd and 100th finite element. 

Similar to the previous example, the bi-active multipliers $v^{REG}$ at the two finite elements (Step 2 of Algorithm) are both appropriately scaled and bounded with ($v^{max} = 161.78$), which ensures that the solution is S-stationary, and thus B-stationary. 


\section{Conclusions and Future Work}
In this study we present a methodology to solve the hybrid dynamics optimal control problem as a dynamic complementarity system (DCS) problem using complementarity constraints, and incorporating moving finite elements with switch detection constraints. The resulting large scale complementarity constraint problem is solved using the active-set strategy to ensure global convergence to local optimal (i.e., B-stationary) solutions. Hybrid dynamic optimization examples show that the proposed method is able to locate accurately the non-smoothness of solutions, even with sliding mode characteristics. In future work, we will use our algorithm to solve different types of optimal control problems with hybrid dynamics, such as in the NOSNOC \cite{nurkanovic2022nosnoc} benchmark.

\bibliography{references}

\begin{thebibliography}{34}
\providecommand{\natexlab}[1]{#1}
\providecommand{\url}[1]{\texttt{#1}}
\expandafter\ifx\csname urlstyle\endcsname\relax
  \providecommand{\doi}[1]{doi: #1}\else
  \providecommand{\doi}{doi: \begingroup \urlstyle{rm}\Url}\fi

\bibitem[Baumrucker and Biegler(2009)]{Baumrucker2009mpec}
B.~Baumrucker and L.~T. Biegler.
\newblock {MPEC} strategies for optimization of a class of hybrid dynamic systems.
\newblock \emph{Journal of Process Control}, 19\penalty0 (8):\penalty0 1248--1256, 2009.

\bibitem[Bemporad and Morari(1999)]{bemporad1999control}
A.~Bemporad and M.~Morari.
\newblock Control of systems integrating logic, dynamics, and constraints.
\newblock \emph{Automatica}, 35\penalty0 (3):\penalty0 407--427, 1999.

\bibitem[Bengea and DeCarlo(2005)]{bengea2005optimal}
S.~C. Bengea and R.~A. DeCarlo.
\newblock Optimal control of switching systems.
\newblock \emph{Automatica}, 41\penalty0 (1):\penalty0 11--27, 2005.

\bibitem[Bertsekas and Ozdaglar(2002)]{bertsekas-ozdaglar}
D.~P. Bertsekas and A.~E. Ozdaglar.
\newblock Pseudonormality and a {Lagrange} multiplier theory for constrained optimization.
\newblock \emph{JOTA}, 114:\penalty0 287--343, 2002.

\bibitem[Bi{\'a}k et~al.(2013)Bi{\'a}k, Hanus, and Janovsk{\'a}]{biak2013some}
M.~Bi{\'a}k, T.~Hanus, and D.~Janovsk{\'a}.
\newblock Some applications of filippov’s dynamical systems.
\newblock \emph{J. Comp. and App. Math.}, 254:\penalty0 132--143, 2013.

\bibitem[Buss et~al.(2002)Buss, Glocker, Hardt, Von~Stryk, Bulirsch, and Schmidt]{buss2002nonlinear}
M.~Buss, M.~Glocker, M.~Hardt, O.~Von~Stryk, R.~Bulirsch, and G.~Schmidt.
\newblock Nonlinear hybrid dynamical systems: modeling, optimal control, and applications.
\newblock In \emph{Modelling, Analysis, and Design of Hybrid Systems}, pages 311--335. Springer, 2002.

\bibitem[Cassandras et~al.(2001)Cassandras, Pepyne, and Wardi]{cassandras2001optimal}
C.~G. Cassandras, D.~L. Pepyne, and Y.~Wardi.
\newblock Optimal control of a class of hybrid systems.
\newblock \emph{IEEE Trans. Aut. Cont.}, 46\penalty0 (3):\penalty0 398--415, 2001.

\bibitem[Filippov(1960)]{filippov1960differential}
A.~F. Filippov.
\newblock Differential equations with discontinuous right-hand side.
\newblock \emph{Matematicheskii sbornik}, 93\penalty0 (1):\penalty0 99--128, 1960.

\bibitem[Flegel and Kanzow(2002)]{flegel-kanzowacq}
M.~L Flegel and C.~Kanzow.
\newblock An {Abadie}-type constraint qualification for math programs with equilibrium constraints.
\newblock Technical report, Inst. Applied Mathematics and Statistics, Univ. W$\ddot u$rzburg, 2002.

\bibitem[Fletcher et~al.(2006)Fletcher, Leyffer, Ralph, and Scholtes]{fletcher}
R.~Fletcher, S.~Leyffer, D.~Ralph, and S.~Scholtes.
\newblock Local convergence of {SQP} methods for mathematical programs with equilibrium constraints.
\newblock \emph{SIAM J. Opt.}, 17\penalty0 (1):\penalty0 259--286, 2006.

\bibitem[Fukushima and Pang(1999)]{fukushima-pang}
M.~Fukushima and J-S. Pang.
\newblock Convergence of a smoothing continuation method for mathematical progams with complementarity constraints.
\newblock In \emph{Ill-posed Variational Problems and Regularization Techniques}. Springer-Verlag, 1999.

\bibitem[Guo et~al.(2015)Guo, Lin, and Ye]{guo-lin-ye}
L.~Guo, G-H. Lin, and J.~Ye.
\newblock Solving mathematical programs with equilibrium constraints.
\newblock \emph{JOTA}, 166:\penalty0 234--256, 2015.

\bibitem[Hedlund and Rantzer(1999)]{hedlund1999optimal}
S.~Hedlund and A.~Rantzer.
\newblock Optimal control of hybrid systems.
\newblock In \emph{Proc. 38th IEEE CDC (Cat. No. 99CH36304)}, volume~4, pages 3972--3977. IEEE, 1999.

\bibitem[Kazi et~al.(2024)Kazi, Thombre, and Biegler]{kazi-thombre-biegler}
S.~R. Kazi, M.~Thombre, and L.~T. Biegler.
\newblock Globally convergent method for optimal control of hybrid dynamical systems.
\newblock In \emph{12th IFAC Intl. Symp. on Adv. Cont. Chem. Proc.}, 2024.

\bibitem[Leyffer(2000)]{leyffer2000}
S.~Leyffer.
\newblock {MacMPEC}: {AMPL} collection of {MPECs}.
\newblock 2000.

\bibitem[Moudgalya and Jaguste(2001)]{moudgalya2001class-II}
K.~M. Moudgalya and S.~Jaguste.
\newblock A class of discontinuous dynamical systems {II}. {An} industrial slurry high density polyethylene reactor.
\newblock \emph{Chem. Eng. Sci.}, 56\penalty0 (11):\penalty0 3611--3621, 2001.

\bibitem[Moudgalya and Ryali(2001)]{moudgalya2001class-I}
K.~M. Moudgalya and V.~Ryali.
\newblock A class of discontinuous dynamical systems {I}. {An} ideal gas--liquid system.
\newblock \emph{Chem. Eng. Sci.}, 56\penalty0 (11):\penalty0 3595--3609, 2001.

\bibitem[Nurkanovi{\'c} and Diehl(2022)]{nurkanovic2022nosnoc}
A.~Nurkanovi{\'c} and M.~Diehl.
\newblock {NOSNOC}: A software package for numerical optimal control of nonsmooth systems.
\newblock \emph{IEEE Control Systems Letters}, 6:\penalty0 3110--3115, 2022.

\bibitem[Nurkanovi{\'c} and Leyffer(2025)]{nurkanovic2025hybrid}
A.~Nurkanovi{\'c} and S.~Leyffer.
\newblock A globally convergent method for computing {B}-stationary points of mathematical programs with equilibrium constraints.
\newblock \emph{arXiv preprint arXiv:2501.13835}, 2025.

\bibitem[Nurkanovi{\'c} et~al.(2024)Nurkanovi{\'c}, Sperl, Albrecht, and Diehl]{nurkanovic2022finite}
A.~Nurkanovi{\'c}, M.~Sperl, S.~Albrecht, and M.~Diehl.
\newblock Finite elements with switch detection for direct optimal control of nonsmooth systems.
\newblock \emph{Numerische Mathematik}, 156:\penalty0 1115--1162, 2024.

\bibitem[Patel et~al.(2019)Patel, Shield, Kazi, Johnson, and Biegler]{patel2019contact}
A.~Patel, S.~Shield, S.~Kazi, A.~M. Johnson, and L.~T. Biegler.
\newblock Contact-implicit trajectory optimization using orthogonal collocation.
\newblock \emph{IEEE Robotics and Automation Letters}, 4\penalty0 (2):\penalty0 2242--2249, 2019.

\bibitem[Raghunathan et~al.(2004)Raghunathan, Diaz, and Biegler]{raghunathan2004MPEC}
A.~U. Raghunathan, M.~S. Diaz, and L.~T. Biegler.
\newblock An {MPEC} formulation for dynamic optimization of distillation operations.
\newblock \emph{Comp. Chem. Eng.}, 28\penalty0 (10):\penalty0 2037--2052, 2004.

\bibitem[Raghunathan et~al.(2022)Raghunathan, Jha, and Romeres]{raghunathan2022pyrobocop}
A.~U. Raghunathan, D.~K. Jha, and D.~Romeres.
\newblock Pyrobocop: Python-based robotic control \& optimization package for manipulation.
\newblock In \emph{2022 Intl. Conf. Robotics Aut. (ICRA)}, pages 985--991. IEEE, 2022.

\bibitem[Riedinger et~al.(1999)Riedinger, Zanne, and Kratz]{Riedinger1999maximum}
P.~Riedinger, C.~Zanne, and F.~Kratz.
\newblock Time optimal control of hybrid systems.
\newblock In \emph{Proc. 1999 Am. Control Conference (Cat. No. 99CH36251)}, volume~4, pages 2466--2470 vol.4, 1999.
\newblock \doi{10.1109/ACC.1999.786491}.

\bibitem[Riedinger et~al.(2003)Riedinger, Iung, and Kratz]{riedinger2003optimal}
P.~Riedinger, C.~Iung, and F.~Kratz.
\newblock An optimal control approach for hybrid systems.
\newblock \emph{European J. Control}, 9\penalty0 (5):\penalty0 449--458, 2003.

\bibitem[Scheel and Scholtes(2000)]{scheel2000mathematical}
H.~Scheel and S.~Scholtes.
\newblock Mathematical programs with complementarity constraints: Stationarity, optimality, and sensitivity.
\newblock \emph{Math. of Op. Res.}, 25\penalty0 (1):\penalty0 1--22, 2000.

\bibitem[Scholtes(2001)]{scholtes2001}
S.~Scholtes.
\newblock Convergence properties of a regularization scheme for math programs with complementarity constraints.
\newblock \emph{SIAM J. Opt.}, 11\penalty0 (4):\penalty0 918--936, 2001.

\bibitem[Shikhman and L\"ammel(2025)]{laemmel-shikhman}
V.~Shikhman and S.~L\"ammel.
\newblock Anomalies of the scholtes regularization for mathematical programs with complementarity constraints.
\newblock \emph{arXiv preprint arXiv:2501.07383v1}, 2025.

\bibitem[Stewart(1990)]{stewart1990high}
D.~Stewart.
\newblock A high accuracy method for solving odes with discontinuous right-hand side.
\newblock \emph{Numerische Mathematik}, 58:\penalty0 299--328, 1990.

\bibitem[Stewart and Anitescu(2010)]{stewart2010optimal}
D.~E. Stewart and M.~Anitescu.
\newblock Optimal control of systems with discontinuous differential equations.
\newblock \emph{Numerische Mathematik}, 114:\penalty0 653--695, 2010.

\bibitem[Van Der~Schaft and Schumacher(2007)]{van2007introduction}
A.J Van Der~Schaft and H.~Schumacher.
\newblock \emph{An introduction to hybrid dynamical systems}, volume 251.
\newblock springer, 2007.

\bibitem[Wang and Biegler(2025)]{wang-biegleronline}
K.~Wang and L.~T Biegler.
\newblock Piecewise {M}-stationarity of local minimizers of {MPCCs} and convergence of {NCP}-based bounding methods.
\newblock \emph{Optimization Online}, https://optimization-online.org/?p=24543, 2025.

\bibitem[Wei et~al.(2008)Wei, Zefran, and DeCarlo]{wei2008optimal}
S.~Wei, M.~Zefran, and R.~A. DeCarlo.
\newblock Optimal control of robotic systems with logical constraints: Application to uav path planning.
\newblock In \emph{2008 IEEE Intl. Conf. Robotics Aut.}, pages 176--181. IEEE, 2008.

\bibitem[Xu and Antsaklis(2004)]{xu2004optimal}
X.~Xu and P.~J Antsaklis.
\newblock Optimal control of switched systems based on parameterization of the switching instants.
\newblock \emph{IEEE Trans. Aut. Cont.}, 49\penalty0 (1):\penalty0 2--16, 2004.

\end{thebibliography}

\end{document}